\documentclass[10pt,a4paper]{article}
\usepackage{graphicx}
\usepackage{amsmath}
\usepackage{amssymb}
\usepackage{mathrsfs}
\usepackage{hyperref}
\usepackage{mathtools}
\usepackage[ruled]{algorithm2e}
\usepackage{cite}
\date{\today}

\usepackage{geometry}
\allowdisplaybreaks[4]

\begin{document}
\title{\bf Deep BSDE-ML Learning and Its Application to Model-Free Optimal Control\thanks{This work is supported in part by the National Natural Science Foundation of China (61773222, 11871369, 61973172, 62173191). }}

\author{Yutian Wang\thanks{College of Artificial Intelligence, Nankai University, Tianjin 300350, P.R. China. Email: {\tt wangyt2239 @mail.nankai.edu.cn}. }~~~~~~~  Yuan-Hua Ni\thanks{College of Artificial Intelligence, Nankai University, Tianjin 300350, P.R. China. Email: {\tt yhni@nankai.edu.cn}. }}
\maketitle

\begin{abstract}

A modified Deep BSDE (backward differential equation) learning method with \emph{measurability loss}, called Deep BSDE-ML method, is introduced in this paper to solve a kind of linear decoupled forward-backward stochastic differential equations (FBSDEs), which is encountered in the policy evaluation of learning the optimal feedback policies of a class of stochastic control problems. The measurability loss is characterized via the measurability of BSDE's state at the forward initial time, which differs from that related to terminal state of the known Deep BSDE method.
Though the minima of the two loss functions are shown to be equal, this measurability loss is proved to be equal to the expected mean squared error between the true diffusion term of BSDE and its approximation. This crucial observation extends the application of the Deep BSDE method---approximating the gradients of the solution of a partial differential equation (PDE) instead of the solution itself.

Simultaneously, a learning-based framework is introduced to search an optimal feedback control of a deterministic nonlinear system. Specifically, by introducing Gaussian exploration noise, we are aiming to learn a robust optimal controller under this stochastic case. This reformulation sacrifices the optimality to some extent, but as suggested in reinforcement learning (RL) exploration noise is essential to enable the model-free learning. The new stochastic optimal control problem is solved with general policy iteration methodology---repeating policy evaluation and policy improvment. Instead of fitting the value function, our policy evaluation approximates its gradient, thus can be seamlessly integrated with policy improvement without manually differentiating a neural network. By using the proposed Deep BSDE-ML method, this is achieved through optimizing the understood loss function in the FBSDE formulation. With some simulating tricks, the whole algorithm can be implemented in both model-based and model-free fashions. Compared with the Markov framework of RL, our method is built on the diffusion process, thus is preferred from a theoretical point of view for continuous-time and continuous-space tasks. Numerical Experiments suggest that the proposed model-free approach performs as good as its model-based counterpart.
\end{abstract}


\section{Introduction}
\label{sec:orgf89f085}

Dynamic programming principle is a fundamental and powerful tool for solving optimal control problems; its basic idea is to consider a family of problems with different initial times and states, and to establish relationships among these problems. Under some mild conditions, the optimal value function satisfies the so-called Hamilton-Jacobi-Bellman (HJB) equation. In most situations, the analytical solutions to HJB equations are hardly obtained due to the equation's high nonlinearity. Hence, numerical solution is much desirable to such kind of nonlinear PDEs. Yet, the classical mesh-based algorithms are computationally expensive due the known ``curse of dimensionality", namely, the computational complexity will exponentially explode with regard to the underlying dimension. Sadly, when the the dimension is bigger than 3, traditional mesh-based algorithms are impractical to some extent.

Recently, deep neural networks have been applied to develop numerical methods that have demonstrated remarkable performance in overcoming the curse of dimensionality and solving high-dimensional HJB equation effectively \cite{E2020,Nakamura,Sirignano2018,Han2018,Onken2021,Pereira2019}. \cite{Nakamura} proposes to use traditional numerical methods to evaluate the solution of HJB equation at certain points, and then use these pre-computed data to train a neural network; relying on the generalization of neural networks, the numerical solution on the entire domain is then obtained. Although being very intuitive, their methods requires enough training data to prevent the overfitting phenomenon, which is expensive in most cases. Deep Galerkin Method (DGM) is claimed to be a general approach to solving PDEs numerically \cite{Sirignano2018}. DGM regards solving a PDE as a variational optimization problem, where the optimized objective function is composed of the PDE itself and boundary conditions. By parameterizing a trial solution of the PDE with a neural network, some gradient optimization method is used to update the neural network so that the deviation from the differential eqaution and boundary condition is reduced. \cite{Han2018} propose the Deep BSDE method to solve a class of parabolic second-order PDEs (HJB equation falls into this category). By applying nonlinear Feynman-Kac's lemma, the problem is transformed into solving a FBSDE, which is further reformulated to a constrained variational optimization problem. Neural networks and gradient descent methods are employed to do the final optimization.

The DGM and Deep BSDE method are two parallel numerical frameworks for solving general PDEs (not limited to HJB equation). These methods do not rely on the training data provided by some external algorithms, thus can be considered as unsupervised learning methods. One drawback of DGM is the high computational complexity of its loss function, which includes the differential operator in the PDE to be solved. On the other hand, the loss function used by Deep BSDE involves only simple additive and other function calculations, thereby Deep BSDE iterates faster. However, because the nonlinear Feynman-Kac's lemma is involved in the reformulation procedure, Deep BSDE method can only handle some specific PDEs, while DGM is a more general approach. Some other works, e.g., \cite{Onken2021} also parameterize the solution of HJB equation with a neural network, but the objective function to be optimized is directly chosen as the cost functional plus some regularization term. And the regularization term is then selected to be the deviation of the trial solution from the PDE and boundary conditions, which is in fact the objective function in DGM. \cite{Pereira2019} solves several HJB equations of a class of nonlinear control systems with Deep BSDE method, but the form of FBSDEs has been modified so that the solution of the forward stochastic differential equation (FSDE) is no longer a trajectory driven by pure noise. To some extent, \cite{Onken2021} and \cite{Pereira2019} can be viewed as the corresponding advanced works of DGM and Deep BSDE in optimal control problems. They no longer consider solving a general PDE, but combine some other properties of HJB eqaution that are inherited from optimal control theory.

Finding the optimal control via solving HJB equation is generally regarded as an indirect approach, because it does not optimize the cost functional directly. Compared with optimization-based methods that find a series of suboptimal policies with decreasing cost, solving HJB equation can give the optimal state-feedback policy once and for all. However, this indirect approach has some drawbacks in practice. In some situations, the optimal policy might be not applicable due to physical constraints or safety concerns. Strictly speaking, for constrained optimal control problems, there are cases that the policy having minimal cost might not satisfy the constraint and thus is excluded from consideration. In such situations, solving HJB equation is not a good choice. Another practical issue is the numerical error. Despite the difficulty of solving general nonlinear PDEs, the obtained policy by solving HJB equation certainly deviates from the optimal one. Therefore, there is no guarantee that this deviated policy works well, especially in nonlinear systems. Therefore, directly applying numerical methods on HJB equation may not yield a satisfactory result.

Inspired by the policy iteration algorithm in reinforcement learning for Markov decision process (MDP) and recent theoretical works for diffusion process, the (nonlinear) HJB equation will be approximated in this paper with a series of linear PDEs, each of which is restated with a linear FBSDE \cite{RL2018,Lee2021,Jia2021SSRN,Wang2020JMLR,Wang2020MAFI}. We regard this approach as the policy iteration in countinous time and space. Under this formulation, policy evaluation is equivalent to solving a linear decoupled FBSDE. Recently, \cite{Jia2021SSRN} proposes a martingale approach to approximate a part of the solution to linear decoupled FBSDE, and arrives at a continuous version of gradient Monte Carlo algorithm. Interestingly, the Deep BSDE method can be shown to approximate the another part of the solution to linear FBSDEs. Based on this observation, we develop a modified Deep BSDE method with \emph{measurability loss}, called Deep BSDE-ML method, to match the model-free martingale approach. Indeed, all the methods mentioned above except the martingale approach are model-based, and requires an exact description of the system dynamics, which is unrealistic in many practical applications. The numerical framework proposed in this paper can be model-free and works well with any black-box simulating environments. This improvement is essential because it makes that those control problems with unknown or complex system can be solved in a clean and reliable way.

Compared with the original Deep BSDE method, the proposed Deep BSDE-ML method utilizes the adaptiveness of the solution, whereas the former concerns the terminal condition. In practice, in order to apply the Deep BSDE method, one starts with an initial guess and an approximated diffusion term. After forward integration, one can adjust the parameters of neural network by checking the deviation from the true terminal condition. The methodology is same for Deep BSDE-ML except for the backward integration and a different loss function, which is measured via the state's measurability at the initial time. We rigorously prove that these two methods are equivalent for linear BSDEs. Interesting, the proposed measurability loss is shown to be strictly equal to the expected mean squared error between the true diffusion term and its approximation. To the best of our knowledge, we are the first to give such an analytical interpretation of the Deep BSDE method for linear BSDEs.

The main contribution of this paper is the introduction of Deep BSDE-ML method and its application to learning an optimal control of a model-free deterministic optimal control problem. Concerned with this control problem, some Gaussian noise is firstly employed to perturb the system and the original optimal control problem is approximated by another stochastic optimal control problem. Following the general policy iteration methodology of RL, we repeatedly perform policy evaluation and policy improvement. In particular, the Deep BSEE-ML method is developed so as to evaluate the policy performance, and both a model-based and a model-free schemes are developed to solve the proposed optimization problem. Compared with the aforementioned works, our contributions are threefold.

\begin{enumerate}
\item The proposed measurability loss is proved to be strictly equal to the expected mean squared error between the diffusion term in BSDE and its approximation. This analytical result puts some insights of the Deep BSDE method for linear FBSDEs.
\item A carefully designed model-free implementation is presented. By introducing some exploration noise, the transformed stochastic problem makes it possible to be implemented totally model free. The adding stochastic noises also makes the obtained policy more robust.
\item On the evaluating the gradient of value function, existing reinforcement learning methods all approximate the value function firstly, and unintentionally differentiate it to obatin the policy gradient. However, our policy evaluation tries to parameterize the gradient of value function and optimize the approximation error directly. This makes the proposed method and existing works are essentially different.
\end{enumerate}

In short, out method works in a policy iteration fashion, where the policy evaluation is performed via the Deep BSDE-ML method and policy improvement is guided by the Hamiltonian. The rest of this paper is organized as follows. In Section 2, we gives the mathematical statement of the considered nonlinear optimal control problem and some preliminaries about HJB equation. In Section 3, the Deep BSDE-ML method is presented together some theoretical justifications and its role in recent advances. Section 4 poses the theory foundation of policy improvement. The model-based and model-based numerical schemes are discussed together some numerical examples in Section 5. Finally, we concludes in Section 6 with some future directions.

\section{Problem formulation and preliminaries}
\label{sec:org2108542}

The optimal control problem considered in this paper is mathematically described as
\begin{align}
\min_{u(\cdot)}\,~&J(u(\cdot)) \coloneqq
\phi(x(T)) + \int_0^T \bigl[Q(t, x(t)) + \frac{1}{2}u^\intercal(t) R(t,x(t)) u(t)\bigr]\,dt,
\label{eq:op}\\
\operatorname{s.t.}\,~& \dot{x}(t) =
F(t, x(t)) + G(t, x(t)) u(t),
\quad x(0) = x_0,~t\in[0,T],
\label{eq:dyn}
\end{align}
where \(T < \infty\) is the time horizon, and \(x(t) \in \mathbb{R}^{D_x}\), \(u(t) \in \mathbb{R}^{D_u}\)  are the state and control input, respectively, with \(x_0\) a fixed initial state. \(J(u(\cdot))\) is the cost functional, composed of a terminal cost \(\phi(x(T))\) and an integral of the running cost. The running cost contains a quadratic term \(u^T(t)R(t,x(t))u(t)/2\) with \(R(t,x(t))\) a positive definite matrix depending on time and state, and a nonlinear term \(Q(t,x(t))\) that is independent of control input. Eq. \eqref{eq:dyn} is the deterministic system dynamic equation, which does not contain any noise and the control dependency is described by a time-variant \(D_x\times D_u\)-valued function \(G(t, x(t))\). For notation simplicity, we sometimes drop the dependency of these functions when it does not cause ambiguity.

The optimal control problem \eqref{eq:op}-\eqref{eq:dyn} can be viewed as a nonlinear extension of the known  linear-quadratic (LQ) problems. The quadratic form of control is retained in cost functional \eqref{eq:op}, but the dependency of states is formulated as some general nonlinear function \(Q(t, x(t))\). The methodology developed in this paper relies on this special model structure, which makes it easier to infer the model information by imposing disturbances on controls. This affine-control structure also makes it possible to find the explicit minimizer of the Hamiltonian. Nevertheless, a large class of control problems in practical applications falls within this category.

The term ``model-free'' mentioned above means that the dynamic equation \eqref{eq:dyn} is unknown to us, namely, the function \(F(t,x), G(t,x)\) is not known. Yet we assume that the cost functional \eqref{eq:op} are given as it is often determined by users. This is the common case in practical engineering applications, e.g., if some key parameters of the system cannot be measured in advance, then the accurate model remains unknown. Instead of requiring an exact mathematical model description, we assume that there exists some black-box environment to be interacted and observed. This environment takes control input and returns current state information. In other words, although we have no idea about how this system involves, we can observe the resulted trajectories. In practice, such a black-box environment appears as a simulating software or a real-world measurement device attached to the concerned system.

However, we shall point out that a deterministic environment is not very suitable for model-free learning, and some stochastic noise for exploration is inevitable. Indeed, when facing a unknown environment, it is better to try some different policies simultaneously, i.e., mixed policies or called stochastic policies, instead of sticking to a deterministic policy. This is the intuitive explanation of the concept of exploration in reinforcement learning: after trying different actions, the benefits caused by these actions can be compared, and consequently the agent (controller) knows better about how to improve his policy. We employ this idea of exploration with an additional uncertainty term in a system dynamics \eqref{eq:dyn}. To be concrete, instead of seeking the optimal controller of the original deterministic problem \eqref{eq:op}, we aim to learn an optimal state-feedback controller under an artificial stochastic system
\begin{equation}
dx(t) = (F(t,x(t))+G(t,x(t))u(t, x(t)))\,dt + \sigma(t,x(t))\,dW(t),\qquad x(0)=x_0,
\label{eq:sto-dyn}
\end{equation}
on a filtered probability space \((\Omega,\mathcal{F},\{\mathcal{F}_t\}_{0\leq t\leq T}, \mathbb{P})\) satisfying the usual condition. \(W:[0,T]\times\Omega\to\mathbb{R}^{D_w}\) the standard Brownian motion and \(\{\mathcal{F}_t, 0\leq t\leq T\}\) is the natural filtration generated by it with \(\mathcal{F}_0\) contains all \(\mathbb{P}\)-null sets. The diffusion coefficient \(\sigma(t,x)\) is some \(D_x\times D_w\) matrix value function and is independent to control. Note that to avoid the subtle concept of measurability of the control \(u(t)\) with respect to \(\mathcal{F}_t\), we turn to find the state-feedback form \(u(t, x(t))\). A state-feedback form is more applicable in practice because it is more robust against uncertainties in real applications. Also, the existence of the optimal state-feedback controller is ensured by the dynamical programming theory under some mild conditions.

One may ask how to choose the diffusion coefficient \(\sigma(t,x)\) and how to implement such the stochastic dynamic \eqref{eq:sto-dyn} with a black-box environment of the deterministic version \eqref{eq:dyn}. Generally speaking, if the environment is presented as a open simulating module, then one can choose arbitrary \(\sigma(t,x)\) and obtain the desired trajectory by manually adding \(\sigma(t,x) dW(t)\) to state \(x(t)\) after each discrete time step. However, for a real-world measurement device or a closed simulating environment, setting the system state manually may be impossible. In this case, the uncertainty term \(\sigma(t,x) dW(t)\) is added via a Gaussian distributed noised added to the control output, and the diffusion coefficient \(\sigma(t,x)\) inherits some model-specific term, see the discussion in Section 4 for more details.

Let us consider the stochastic optimal control problem defined by cost functinal \(\mathbb{E}[J(u(\cdot))]\) and the state dynamics \eqref{eq:sto-dyn}.  In pratice, we may assume we search the optimal policy in a parameterized function space \(\{u^\psi(t,x), \psi\in\Psi\}\), termed policy space, where \(\Psi\) is an abstract index set and \(u^\psi\) denote an entity in it. Assume that \(\Psi\) and \(\psi\mapsto u^\psi\) are chosen such that there exists \(\psi^*\) and \(u^{\psi^*}\) is an optimal state-feedback control policy, abbreviated as \(u^*\). Our stochastic optimal control problem then becomes
\begin{align}
\min_{\psi\in\Psi}\quad& \operatorname{\mathbb{E}}\left[\phi(x(T)) + \int_0^T\bigl[Q(t,x(t)) + \frac{1}{2}\langle Ru^\psi(t,x(t)),u^\psi(t,x(t))\rangle\bigr]\,dt\right],\\
\mathrm{s.t.}\quad & dx(t) = (F(t,x(t))+G(t,x(t))u^\psi(t,x(t)))\,dt + \sigma(t,x(t))\,dW_t.
\end{align}
According to the dynamical programming theory, the optimal policy satisfies the following minimum condition \cite{Yong1999}
\begin{equation}
H(t,x(t),u^*(t,x(t)), v^*_x(t,x(t))) = \min_{u}H(t, x(t), u, v^*_x(t,x(t))),
\label{eq:u-star}
\end{equation}
where \(H(t,x,u,p)\) is the Hamiltonian associated with the stochastic optimal control problem
\begin{equation}
H(t,x,u,p) \coloneqq Q(t,x) + \frac{1}{2}u^\intercal R(t,x)u + \langle p, F(t,x) + G(t,x)u\rangle.
\label{eq:def-Hamiltonian}
\end{equation}
And \(v^*_x\) is the gradient of the optimal value function \(v^*(t,x)\) defined as
\begin{equation}
v^*(t,x) \coloneqq \min_{u(\cdot)}\,\mathbb{E}\left[\phi(X_T) + \int_t^TQ(\tau, X_\tau) + \frac{1}{2}u_\tau^\intercal R(\tau,X_\tau) u_\tau\,d\tau\biggm|X_t=x\right],
\label{eq:def-value-function}
\end{equation}
where \(\{X_\tau, t\leq\tau\leq T\}\) is the adapted solution to the stochastic differential equaiton \eqref{eq:sto-dyn} starting from \(X_t=x\). Noticing that the minimum in Eq. \eqref{eq:def-value-function} is achieved at \(u=u^*(t,x)\) and applying the Feynman-Kac's lemma, the above conditional expectation is related to the solution to the following partial differential equation \cite{Karatzas2012}
\begin{equation}
\mathcal{L}^*v(t,x) + Q(t,x) + \frac{1}{2}u^*(t,x)^\intercal Ru^*(t,x) = 0,\qquad v(T,x)=\phi(x),
\label{eq:PDE-formulation}
\end{equation}
where \(\mathcal{L}\)\textsuperscript{*} denotes the infinitesimal generator of \(u^*\)
\begin{equation}
\mathcal{L}^*v(t,x) \coloneqq v_t(t,x) + \langle v_x(t,x), F(t,x)+G(t,x)u^*(t,x)\rangle + \frac{1}{2}\operatorname{tr}\left\{v_{xx}(t,x)\sigma\sigma^\intercal(t,x)\right\}.
\label{eq:def-L-star}
\end{equation}

Hence, Eq. \eqref{eq:u-star} along with \eqref{eq:def-Hamiltonian} has given the explicit form of the optimal state-feedback controller with the help of the optimal value function \(v^*(t,x)\), or more exactly, with its gradient \(v^*_x(t,x)\). This is essentially different from the MDP concerned in reinforcement learning with discrete time and space, where only \(v^*(t,x)\) is evaluated, not its gradient. However, the so called state-action value function in RL does have some interpretation like \(v^*_x\), and approves a policy improvement direction. Indeed, one may show that there is some subtle relationship between the optimal state-action value function and the Hamiltonian as the discrete time step approaches zero.

To find \(v^*_x\), one may seek the solution to PDE \eqref{eq:PDE-formulation}. Sadly, it is not a simple task, as Eq. \eqref{eq:PDE-formulation} is a nonlinear PDE in fact. The nonlinear term is the quadratic form of \(u^*(t,x)\), which is related to \(v_x^*(t,x)\) regarding Eq. \eqref{eq:u-star}. Actually, by substituting Eq. \eqref{eq:u-star}, \eqref{eq:def-Hamiltonian} and \eqref{eq:def-L-star} into Eq. \eqref{eq:PDE-formulation}, one may recover the well known HJB equation
\begin{equation}
v_t + \min_{u}H(t,x,u,v_x^*) + \frac{1}{2}\operatorname{tr}\left\{v_{xx}\sigma\sigma^\intercal(t,x)\right\} = 0,\qquad v(T,x)=\phi(x).
\label{eq:hjb}
\end{equation}
As explained before, solving high-dimensional HJB equation directly is quite difficult, and even impossible in model-free settings. To overcome this difficulty, we need transform it to solving FBSDEs, which admits the possibility of model-free learning.

\section{Policy evaluation---Deep BSDE-ML method}
\label{sec:orga455d0a}

In this section, we initialize the policy iteration in our framework with policy evaluation, and next section presents the policy improvement. As argued in the last section, the HJB equation \eqref{eq:hjb}, whose solution is the optimal value function, is nonlinear and quite difficult to solve either analytically or numerically. However, if the control in Eq. \eqref{eq:def-value-function} is fixed to a known feedback policy \(u^\psi(t,x)\), then the corresponding value function \(v^\psi(t,x)\) is no longer optimal. In exchange, this value function satisfies a linear PDE
\begin{equation}
\mathcal{L}^\psi v^\psi + Q(t,x) + \frac{1}{2}\langle Ru^\psi(t,x), u^\psi(t,x)\rangle,\qquad v^\psi(T,x)=\phi(x),
\label{eq:PDE-formulation-psi}
\end{equation}
where \(\mathcal{L}^\psi\) is defined by replacing the \(u^*\) with \(u^\psi\) in Eq. \eqref{eq:def-L-star}. This PDE is linear because \(u^\psi\) is known as a prior, while the HJB equation \eqref{eq:PDE-formulation} is nonlinear because \(u^*\) is unknown and must be replaced by the quadratic form of the gradient of optimal value function. Now we are ready to derive the FBSDE representation of the value function \(v^\psi\). Let \(\{X_t, 0\leq t\leq T\}\) be the adapted solution to SDE \eqref{eq:sto-dyn} starting at \(X_0=x_0\) with \(u=u^\psi(t,x)\). Then by applying Itô's lemma to differentiate \(v^\psi(t,X_t)\), one has
\begin{equation}
\label{eq:differential-v}
dv^\psi(t,X_t) = \mathcal{L}^\psi v^\psi(t, X_t)\,dt + \langle\sigma^\intercal v^\psi_x(t,X_t),dW(t)\rangle.
\end{equation}
Combining the forward dynamics \eqref{eq:sto-dyn} and eliminating \(\mathcal{L}^\psi v^\psi\) by \eqref{eq:PDE-formulation-psi} yields the following FBSDE
\begin{align}
\label{eq:fsde}
dX_t&=(F(t,X_t)+G(t,X_t)u^\psi(t,X_t))\,dt+\sigma(t,X_t)\,dW_t,&\qquad X_0&=x_0,\\
\label{eq:bsde}
dY_t&=-(Q(t,X_t)+ \frac{1}{2}u^\psi(t,X_t)^\intercal Ru^\psi(t,X_t))\,dt + \langle Z_t, dW_t\rangle,&\qquad Y_T&=\phi(X_T),
\end{align}
which admits the following adapted solution
\begin{equation}
Y_t = v^\psi(t, X_t),\qquad Z_t = \sigma^\intercal v^\psi_x(t, X_t).
\label{eq:FBSDE-formulation}
\end{equation}
Compared with the PDE formulation \eqref{eq:PDE-formulation-psi}, this BSDE characterization is model-free and also derivative-free. The former uses the generator \(\mathcal{L}^\psi\) which involves differentiating operation and model dynamics \(F, G\) too, while the later hides these model-based information in the forward equation and its solution \(X_t\). Hence, we turn our attention from Eq. \eqref{eq:PDE-formulation-psi} to the FBSDE \eqref{eq:fsde}--\eqref{eq:bsde} and are devoted to find an adapted solution. To this aim, a learning scheme is developed in the following subsection, which is termed the Deep BSDE-ML method.

For the sake of simplicity, we set the policy space \(\{u^\psi, \psi\in\Psi\}\) to be the collection of policies with \(C^{1,2}([0,T]\times\mathbb{R}^{D_x})\) value functions, where a function \(v(t,x)\in C^{1,2}\) means \(v_t, v_x, v_{xx}\) are continuous on the interior of its domain. Also, we denote the set of all \(\{\mathcal{F}_t\}_{0\leq t\leq T}\)-adapted process \(\{\xi_t, 0\leq t\leq T\}\) with \(\operatorname{\mathbb{E}}\int_0^T\|\xi_t\|\,dt < \infty\) by \(\mathscr{L}_{\mathcal{F}}^2([0,T])\).

\subsection{Deep BSDE-ML method}
\label{sec:orgfc9cb39}

As it will be discussed soon that our policy improvement needs only \(v_x^\psi\), here the process \(Z_t\) is more interested than \(Y_t\) for our purpose. Let \(z^\theta(t,x)\) be some parameterized function and \(\theta\in\Theta\), where \(\Theta\) is an abstract index set termed value space. Define an auxiliary process
\begin{equation}
y^\theta_t \coloneqq \phi(X_T) + \int_t^Tg^\psi(\tau,X_\tau)\,d\tau - \int_t^T\langle z^\theta(\tau, X_\tau), dW_\tau\rangle,
\label{eq:def-ytheta}
\end{equation}
where we denote the running cost under policy \(u^\psi\) by
\begin{equation}
g^\psi(t,x)\coloneqq Q(t,x)+ \frac{1}{2}\langle Ru^\psi(t,x),u^\psi(t,x)\rangle.
\end{equation}
At the first glance, \((y^\theta_t, z^\theta(t,X_t))\) seems to satisfy the differential equality \eqref{eq:bsde} with terminal condition. Unfortunately, this does not mean that we have found an adapted solution, as \(y^\theta_t\) may be not \(\mathcal{F}_t\)-measurable if one simply integrate backwardly! It is not surprising because \(z^\theta\) is arbitrary and one can not expect it coincidences the true answer. Noting that the adapted solution of the BSDE \eqref{eq:bsde} is a pair of processes \((Y_t,Z_t)\), \(Y_t\) will be measurable to \(\mathcal{F}_t\) if \(Z_t\) is particularly selected. In other words, \(y^\theta_t\) of \eqref{eq:def-ytheta} should be \(\mathcal{F}_t\)-measurable at any time instant \(t\) if it is indeed a part of a true solution. Notice that the expression of \(y_t^\theta\) depends on the future forward state \(\{X_\tau, t < \tau\leq T\}\), which at \(\tau\) is measurable with respect to \(\mathcal{F}_\tau\) (bigger than \(\mathcal{F}_t\)). As a matter of fact, the last term in Eq. \eqref{eq:def-ytheta} is the key to fix this measurability issue. For proper \(z^\theta\), the stochastic integral \(\int \langle z^\theta_t, dW_\tau\rangle\) can eliminate the uncertainties caused by the other terms in the right side of Eq. \eqref{eq:def-ytheta}. This is the reason why a single backward equation should be characterized by a pair of processes.

Conversely, if \(y^\theta_t\) is indeed adapted, then we have found a solution to BSDE \eqref{eq:bsde} and finished the policy evaluation procedure. Recall that the adaptiveness of \(\{y_t^\theta, 0\leq t\leq T\}\) with respect to the filtration \(\{\mathcal{F}_t\}_{0\leq t\leq T}\) means some measurability, which can be characterized by the conditional expectation
\begin{equation}
\mathbb{E}[(y_t^\theta-\mathbb{E}[y_t^\theta\,|\,\mathcal{F}_t])^2\,|\,\mathcal{F}_t]=0,\qquad \forall t\in[0,T].
\label{eq:measurability-yt}
\end{equation}
It is interesting that if Eq. \eqref{eq:measurability-yt} holds at \(t=0\), then it is sufficient enough to claim that \(y_t^\theta\) is adapted, and consequently \eqref{eq:measurability-yt} holds on the whole interval \([0,T]\). One can easily verify this from the following observation. By definition,
\begin{equation}
y_0^\theta = \phi(X_T) + \int_0^Tg^\psi(\tau,X_\tau)\,d\tau - \int_0^T\langle z^\theta(\tau, X_\tau), dW_\tau\rangle.
\label{eq:ytheta-0}
\end{equation}
Subtracting Eq. \eqref{eq:def-ytheta} by \eqref{eq:ytheta-0} leads to
\begin{equation}
y_t^\theta = y_0^\theta - \int_0^tg^\psi(\tau,X_\tau)\,d\tau + \int_0^t\langle z^\theta(\tau, X_\tau), dW_\tau\rangle.
\label{eq:ytheta-forward}
\end{equation}
If Eq. \eqref{eq:measurability-yt} holds at \(t=0\), then \(y_0^\theta\) is \(\mathcal{F}_0\)-measurable, which implies that \(y_t^\theta\) is \(\mathcal{F}_t\)-measurable. This simple observation inspires us to consider the following loss function with respect to \(\theta\), referred as measurability loss
\begin{equation}
\mathit{MeasurLoss}(\theta)\coloneqq \operatorname{\mathbb{E}}|y_0^\theta - \operatorname{\mathbb{E}}y_0^\theta|^2.
\label{eq:def-measurability-loss}
\end{equation}
If there exists \(\theta^*\) such that the measurability loss at that point is exactly zero, then \(y_0^{\theta^*}\) is equal to its expectation almost surely. This implies \(\{y_0^\theta\neq \operatorname{\mathbb{E}}y_0^\theta\}\) is a \(\mathbb{P}\)-null set thus is contained in \(\mathcal{F}_0\). Hence, \(y_0^\theta\) is \(\mathcal{F}_0\)-measurable and \(y_t^\theta\) is \(\mathcal{F}_t\)-measurable by Eq. \eqref{eq:ytheta-forward}. Some careful analysis reveals that this measurability loss quantifies the deviation from the ground true \(Z=\sigma^\intercal v_x^\psi\) along the state trajectory \(\{X_t, 0\leq t\leq T\}\).

\textbf{Theorem 1}.  Fix a policy \(u^\psi\) and denote its value function by \(v^\psi\). Suppose \(v^\psi\in C^{1,2}\) and there exists an unique strong solution to the forward SDE \eqref{eq:fsde}, denoted by \(\{X_t, 0\leq t\leq T\}\). Let \(\{y_t^\theta, 0\leq t\leq T\}\) denote the process constructed with \(z^\theta(\cdot,\cdot)\) by Eq. \eqref{eq:def-ytheta}. If \(\{z^\theta(t,X_t):0\leq t\leq T\}\in\mathscr{L}_{\mathcal{F}}^2([0,T])\), then
\begin{equation}
\operatorname{\mathbb{E}}|y_0^\theta - \operatorname{\mathbb{E}} y_0^\theta|^2 = \operatorname{\mathbb{E}}\int_0^T\|z^\theta(t,X_t) - \sigma^\intercal v_x^\psi(t,X_t)\|^2\,dt.
\label{eq:z-err}
\end{equation}

\textbf{Remark}.  This theoretical result gives a new way to think about the measurability loss. In deriving the measurability loss, we can only say that if this loss decreases to zero, then we find a true solution of the BSDE. However, it is not clear before this theorem that how good the approximation is during the optimization. It may be the case that the measurability loss do have a strictly minimizer at the point we want but behaves totally different otherwhere. Fortunately, Eq. \eqref{eq:z-err} dispels this worries. It points out that the measurability loss is strictly equal to the expected mean squared error between the parameterized function and the ground truth along the state trajectory. This equality also ensures that if the loss is zero, then \(z^\theta(t,X_t)\) equals \(\sigma^\intercal v_x^\psi(t,X_t)\) almost surely. Integrating yields the conclusion that \(y_t^\theta\) equals \(v^\psi(t,X_t)\) almost surely too. However, we shall also point out this theoretical result heavily relies on the fact that the running cost \(g^\psi\) is independent of \(Y_t\) and \(Z_t\). This might not hold for nonlinear BSDEs. Nevertheless, the measurability loss defined by Eq. \eqref{eq:def-ytheta} and \eqref{eq:def-measurability-loss} is a considerably rational choice for approximating the value function gradients.

\subsection{Related to recent advances}
\label{sec:org1c60f0c}

We are not the first to formulate the policy evaluation as solving a decoupled linear FBSDE. \cite{Jia2021SSRN} consider the same FBSDE representation \eqref{eq:fsde}--\eqref{eq:FBSDE-formulation} too. Their goal is approximating the value function itself by some parameterized function \(J^\theta(t,x)\). Despite the different choice of target function, their starting point is same as ours. With a trial solution \(\widetilde{Y}_t = J^\theta(t,X_t)\) at hand, they aim to ensure there exists an adapted process \(\widetilde{Z}_t\) such that \((\widetilde{Y}_t,\widetilde{Z}_t)\) satisfies the BSDE \eqref{eq:bsde}. Suppose such a process indeed exists, and integrating Eq. \eqref{eq:bsde} and rearranging yield
\begin{equation}
\int_0^t\langle\widetilde{Z}_t,dW_t\rangle = J^\theta(t,X_t) + \int_0^tg^\psi(\tau,X_\tau)\,d\tau - J^\theta(0,x_0),\qquad\forall t\in[0,T].
\label{eq:martingale-interpretation}
\end{equation}
By martingale representation theorems, the process \(\widetilde{Z}_t\) exists if and only if the right hand side of Eq. \eqref{eq:martingale-interpretation} is a martingale. Based on this, they propose a loss function, called martingale loss, to ensure this martingality. Then a result similar to Eq. \eqref{eq:z-err} can be established. That is, their martingale loss strictly is strictly equal to the expected mean squared error between the parameterized function \(J^\theta\) and the ground truth \(v^\psi\) along the state trajectory plus a term independent of the choice of \(J^\theta\).

Compared with our result, \cite{Jia2021SSRN} concentrate on finding the true \(Y_t\), while this paper is on finding the true \(Z_t\). This difference comes from different purposes. If one wants to evaluate the value function itself, then \(Y_t\) is a more desirable target. In contrast, we need value function gradients so prefer fitting \(Z_t\) directly. Of course, if one can perfectly solve the optimization problem and decrease the loss function to zero, then both methods find a part of the true solution, and by some differentials or integrals one can represent the other part of the solution. As we have pointed out, as long as Eq. \eqref{eq:z-err} is zero, the auxiliary process \(y_t^\theta\) defined by Eq. \eqref{eq:def-ytheta} turns out to be the true solution \(Y_t\). Similarly, once the martingality condition derived from Eq. \eqref{eq:martingale-interpretation} is satisfied by some \(J^\theta\), then \(\sigma^\intercal J^\theta_x(t,X_t)\) becomes the true solution of \(Z_t\). Regrettably, the exact minimizer of the optimization problem is generally impossible to be found in most practical situations. When approximation error is inevitable, it is preferable to optimize the distance between the parameterized function and the target directly, and avoids the consequent integration or differentiation. Therefore, the methods proposed in \cite{Jia2021SSRN} and this paper are two complementary ways to solve a pair of solution to linear BSDE---one for \(Y_t\) and the other for \(Z_t\).

As discussed at the last of Section 1, our policy evaluation procedure is based on the Deep BSDE method. \cite{Han2018} propose the Deep BSDE method to solve a wide class of second order nonlinear parabolic PDEs. After applying the nonlinear Feynman-Kac's formula, they arrive at a decoupled FBSDE too. Although their ultimate goal is same as \cite{Jia2021SSRN}---trying to find the solution to the PDE---\cite{Han2018} choose to replace the true \(Z_t\) with \(z^\theta(t,X_t)\). This is the reason why we claim our method is based on Deep BSDE. However, our theoretical result suggests that this parameterizing strategy fits the value function gradients, not the value function itself. \cite{Yu2020} has mentioned this variance form loss Eq. \eqref{eq:def-measurability-loss} and its variants, primarily based on considerations of numerical implementations. To the best of our knowledge, we are the first rigorously prove that the Deep BSDE method for linear BSDEs is practically trying to fit the value function gradients instead of the value function itself. The rest of this section briefly review the Deep BSDE method in our settings and prove the equivalence between the this method and ours.

Considering the BSDE \eqref{eq:bsde}, Deep BSDE method wants to evalue the value function at the initial time \(v^\psi(0,x_0)\). In order to achieve this goal, it chooses an initial guess of that value \(y_0^{DB}\) and a parameterized function \(z^\theta(t,x)\). Different from our definition of the \(y_t^\theta\) in Eq. \eqref{eq:def-ytheta}, it tries to integrate Eq. \eqref{eq:bsde} on the time interval \([0,t]\)
\begin{equation}
y_t^{\theta,DB} \coloneqq y_0^{DB} -\int_0^t g^\psi(\tau,X_\tau)\,d\tau + \int_0^tz^\theta(\tau,X_\tau)\,d\tau.
\label{eq:def-ytheta-DB}
\end{equation}
It is clear that this auxiliary process is adapted and satisfies the differential equality in Eq. \eqref{eq:bsde}. Therefore, once \(y_T^{\theta,DB}\) equals the terminal cost \(\phi(X_T)\) almost surely, then one can conclude that the original BSDE has been solved. Based on this necessary condition, Deep BSDE tries to minimize the following loss function
\begin{equation}
\mathit{DBLoss}\coloneqq \operatorname{\mathbb{E}}|y_T^{\theta,DB} - \phi(X_T)|^2
\label{eq:def-DB-loss}
\end{equation}
by adjusting two variables, \(y_0^{DB}\) and \(\theta\). The theoretical analysis in the original paper ends here. However, we can show that for BSDE \eqref{eq:bsde} optimize the loss function Eq. \eqref{eq:def-DB-loss} is equivalent to minimize the measurability loss defined by Eq. \eqref{eq:measurability-yt}. Therefore, Theorem 1 justifies that the Deep BSDE method unexpectedly fits the value function gradients.

\textbf{Theorem 2}. Consider the FBSDE representation Eq. \eqref{eq:fsde}--\eqref{eq:FBSDE-formulation}. Fix a policy \(u^\psi\) and denote its value function by \(v^\psi\). Let the conditions in Theorem 1 holds and \(\{X_t\}\) denotes the strong solution to the FSDE. Let \(y^\theta_t\) and \(y^{\theta,DB}_t\) be two process defined by Eq. \eqref{eq:def-ytheta} and Eq. \eqref{eq:def-ytheta-DB} respectively, where \(\theta\) and \(y_0^{DB}\) are two variables to be optimized. Then it holds that
\begin{equation}
\min_{\theta,y_0^{DB}} \operatorname{\mathbb{E}}|y_T^{\theta,DB} - \phi(X_T)|^2 = \min_{\theta} \operatorname{\mathbb{E}}|y_0^\theta - \operatorname{\mathbb{E}}y_0^\theta|^2.
\end{equation}
Moreover, if the left hand side optimization problem attains its minimum at \((y_0^{*,DB}, \theta^*)\), then the right hand side also attains its minimum at \(\theta^*\). Conversely, if the right hand side attains its minimum at \(\theta^*\), then there exists a \(y_0^{*,DB}\) such that the left hand side optimization problem attains its minimum at \((y_0^{*,DB}, \theta^*)\).

\textbf{Remark}. Compared with the measurability loss used in this paper, Deep BSDE method chooses a loss function depending not only on \(\theta\) but also on another variable \(y_0^{DB}\). We suspect this is because the Deep BSDE method is proposed to find the optimal \(y_0^{*,DB}\), i.e., the initial value \(v^\psi(0,x_0)\), leaving \(\theta\) as an auxiliary variable in optimization. As shown in the proof of above theorem, the loss function of Deep BSDE method can be decomposed into two part, where one part coincides with the measurability loss in this paper and the other part is simply the squared error between \(y_0^{DB}\) and \(v^\psi(0,x_0)\), i.e.,
\begin{equation}
\operatorname{\mathbb{E}}|y_T^{\theta, DB} - \phi(X_T)|^2 = |y_0^{DB} - v^\psi(0,x_0)|^2 + \operatorname{\mathbb{E}}|y_0^\theta - \operatorname{\mathbb{E}}y_0^\theta|^2.
\end{equation}
This decomposition is true for linear BSDEs. Therefore, the loss function concerned by Deep BSDE clearly do fulfill their hope---reflecting the distance between their guess and the value function at the starting point. However, it is seen that Deep BSDE also tries optimizing \(\theta\), or equivalently, optimizing the measurability loss, which is eventually the error of fitting the value function gradients. In general, the Deep BSDE should be recognized as a way which fits the value function gradients and the initial value simultaneously. If there is no need to obtain the initial value, our method is more preferable.

\section{Policy improvement based on value function gradients}
\label{sec:orgb80841e}

Another key of policy iteration method is the policy improvement subroutine. As its name suggests, policy iteration iterates the policy itself. Given an initial policy \(u^\psi\), it manages to find a better policy \(u^{\psi'}\) via the policy improvement procedure. From this perspective, the policy evaluation is only a preliminary procedure to perform policy improvement. One can think the policy iteration as an abstract optimization procedure, where the policy evaluation plays the role of the oracle, which delivers the information about the current policy. Take the gradient-based optimization as an example, which works by repeating evaluating gradient and performing gradient descent. In policy iteration, the methodology is same. Policy evaluation computes some auxiliary information then policy improvement go a step based on it. Usually, it is required that the new policy constructed by policy improvement must be not bad than the old one. In this case, the optimal policy is obviously a fixed point of the policy iteration method.

For deterministic systems, \cite{Lee2021} suggests that a minimization operation over the Hamiltonian yields a policy fulfilling the requirement of policy improvement. In the stochastic settings, \cite{Ni2013} has given a similar result except that the Hamiltonian is replaced by the generalized Hamiltonian. For the problem we are concerned with, minimizing the Hamiltonian is sufficient to obtain a better policy because the diffusion coefficient \(\sigma\) is independent of control. This condition is generally satisfied as \(\sigma\) is actually determined by ourself and it will be discussed soon about how to choose the diffusion coefficient. Following the common procedures to obtain policy improvement, we establish a result for policy comparison first.

\textbf{Theorem 3 }. Given two policy \(u^\psi, u^\theta\), let \(v^\psi,v^\theta\) denote the corresponding value function respectively. Assume \(v^\psi, v^\theta\in C^{1,2}([0,T]\times\mathbb{R}^{D_x})\). Suppose for any \(s\in[0,T]\) and \(x_s\in\mathbb{R}^{D_x}\), there exists a unique strong solution, denoted by \(\{X^\psi_t, s\leq t\leq T\}\), to the state dynamics driven by \(u^\psi\) with initial condition \(x(s)=x_s\). Then
\begin{align}
v^\psi(s,x_s) - v^\theta(s,x_s) &= \operatorname{\mathbb{E}}\int_s^TH(t,X^\psi_t, u^\psi(t,X_t^\psi),v_x^\theta(t,X_t^\psi))\,dt \notag\\
&\qquad\qquad - \operatorname{\mathbb{E}}\int_s^T H(t,X_t^\psi,u^\theta(t,X_t^\psi),v_x^\theta(t,X_t^\psi)\,dt.
\end{align}

This result suggests that for a given policy \(u^\theta\), if there exists a policy \(u^\psi\) such that for any \(t\in[0,T]\) it holds almost surely that
\begin{equation}
H(t,X_t^\psi,u^\psi(t,X_t^\psi),v^\theta_x(t,X_t^\psi)) \leq H(t,X_t^\psi,u^\theta(t,X_t^\psi),v^\theta_x(t,X_t^\psi),
\label{eq:decrease-Hamiltonian}
\end{equation}
then \(u^\psi\) is a policy not bad than \(u^\theta\). Considering the quadratic structure of the Hamiltonian in our settings, we may construct a better policy by setting it to the minimizer of the Hamiltonian. That is, the policy
\begin{equation}
u^\psi(t,x) \coloneqq -R^{-1}G^\intercal v^\theta_x(t,x)
\label{eq:policy-improvement}
\end{equation}
always satisfies the Eq \eqref{eq:decrease-Hamiltonian} and therefore is not bad than \(u^\theta\) according to Theorem 3. This policy improvement strategy has been widely adapted in the theory of adaptive dynamic programming for deterministic systems \cite{Murray2002,ISOCPI2017,Yang2021}. However, Eq. \eqref{eq:policy-improvement} explicitly requires the model-specific information \(G(t,x)\), leaving most of existing works based on it are basically model-based. Some works remove this restriction by training a system identifier simultaneously. In contrast, we accomplish the model-free policy improvement with injecting noise. As suggested in reinforcement learning, exploration noise is the key to be model-free. We sacrifices some optimality in exchange for the free of choosing the diffusion coefficient \(\sigma\).

Now suppose \(\sigma\) is chosen such that their exists matrix valued function \(\Upsilon(t,x)\) satisfying \(G(t,x) = \sigma \Upsilon(t,x)\). Then the policy improvement Eq. \eqref{eq:policy-improvement} can be rewritten as
\begin{equation}
u^\psi(t,x) \coloneqq - R^{-1}\Upsilon^\intercal z^\theta(t,x),\qquad z^\theta(t,x) = \sigma^\intercal v^\theta_x(t,x).
\label{eq:policy-improvement-model-free}
\end{equation}
It is very natural to expect some method to evaluating this \(z^\theta(t,x)\) function. Our policy evaluation is developed precisely to meet this requirement. In this paper, we primarily concern two choices of \(\sigma\) and the corresponding \(\Upsilon\)
\begin{align}
\sigma(t,x) &= \sigma_0I,&\quad&\mathrm{and}&\quad \Upsilon(t,x) &= \sigma_0^{-1}G(t,x),\label{eq:model-based-settings}\\
\sigma(t,x) &= \sigma_0G(t,x),&\quad&\mathrm{and}&\quad \Upsilon(t,x) &= \sigma_0^{-1}I,\label{eq:model-free-settings}
\end{align}
where \(\sigma_0 > 0\) is a constant. In the next section, we will show that the first choice leads to a model-based policy iteration algorithm and the second choice leads to a model-free one.

\section{Numerical Scheme and Examples}
\label{sec:orgbea44a9}

In this section, we first formally describe our policy iteration algorithms and its model-based and model-free implementations, Then we demonstrate how to apply them on a concrete optimal control problems.

Basically, the policy iteration works like this. Firstly, let \(i=0\) and choose an initial policy \(u^{(i)}\), usually a policy that always return zeros. Secondly, we need to evaluate the value function gradient of this policy \(\sigma^\intercal v_x^{(i)}\). Let \((X_t^{(i)},Y_t^{(i)},Z_t^{(i)})\) denote the solution to the FBSDE \eqref{eq:fsde}--\eqref{eq:FBSDE-formulation} where \(u^\psi\) is set to \(u^{(i)}\). In this subroutine, we first numerically solve forward SDE \eqref{eq:FBSDE-formulation}. Then, we fit the diffusion term of the backward SDE \eqref{eq:FBSDE-formulation} with \(z^\theta(t,X^{(i)}_t)\), where \(z^\theta(\cdot, \cdot)\) is a neural network. Thirdly, we construct a new policy \(u^{(i+1)}(t,x)\) based on the minimizer of Hamiltonian, i.e., \(-R^{-1}\Upsilon^\intercal\sigma^\intercal v_x^{(i)}(t,x)\). Repeating updating \(z^\theta\) and policy \(u^{(i)}\) until there is no performance improvement.

To implement this policy iteration procedure with a computer, we need to  discretize the time interval \([0,T]\). Let \(\{t_k:0\leq k\leq H\}\) be uniformly spaced points on \([0,T]\),
\begin{equation}
0 = t_0 < t_1 < t_2 < \cdots < t_{H-1} < t_{H} = T,\qquad t_{k+1}-t_{k} = T/H \eqqcolon \Delta t.
\end{equation}
The forward dynamics driven by a policy \(u^{(i)}(t,x)\) is discretized accordingly,
\begin{equation}
X^{(i)}_{t_{k+1}} - X^{(i)}_{t_k} = (F(t_k,X^{(i)}_{t_k})+G(t_k,X^{(i)}_{t_k})u^{(i)}(t_k,X^{(i)}_{t_k}))\,\Delta t + \sigma(t_k,X^{(i)}_{t_k})\,(W_{t_{k+1}} - W_{t_k}),
\label{eq:model-based-samples}
\end{equation}
where we simulate \(W_{t_{k+1}} - W_{t_k}=\Delta W_k\) by drawing from normal distribution \(\mathcal{N}(0, \Delta t\cdot I)\). However, simulating the state process by Eq. \eqref{eq:model-based-samples} is hard without an exact system model because it requires disturbing the system state with noise \(\sigma \Delta W\) at each time instant. To enable model-free sampling, we employ the disturbances on the control output
\begin{equation}
X^{(i)}_{t_{k+1}}-X^{(i)}_{t_k} = (F(t_k,X^{(i)}_{t_k}) + G(t_k,X^{(i)}_{t_k})(u^{(i)}(t_k,X^{(i)}_{t_k}) + \xi_{t_k}))\Delta t,
\label{eq:model-free-samples}
\end{equation}
where \(\xi_{t_k} = \sigma_0 \Delta W_k/ \Delta t\sim\mathcal{N}(0, \sigma_0^2/\Delta t)\). This case corresponds to the choice of \(\sigma, \Upsilon\) specified in Eq. \eqref{eq:model-free-settings}. Simulating the process determined by Eq. \eqref{eq:model-free-samples} is extremely simple. One need only add a Gaussian distributed noise to current control policy. This can generally be done for both numerical models and real world experiments. From the perspective of reinforcement learning, \(X^{(i)}_{t_k}\) generated via Eq. \eqref{eq:model-free-samples} is nothing but the data collected by an behivour policy, which is Gaussian with mean \(u^{(i)}\) and covariance \(\sigma^2_0/\Delta t\). And our method aims to find the optimal mean \(u^*\) and left the covariance unchanged during training. Notice that if we choose \(\sigma, \Upsilon\) as Eq. \eqref{eq:model-free-settings}, the policy improvement Eq. \eqref{eq:policy-improvement-model-free} becomes model-free too; this results that the whole algorithm does not require the information about model-specific terms \(F(t,x)\) and \(G(t,x)\).

To fit \(z^\theta(t,X^{(i)}_t)\) against \(Z^{(i)}_t\) with sample paths of \(X^{(i)}_t\), we utilize the proposed measurability loss in Eq. \eqref{eq:def-measurability-loss}. Theorem 1 ensures that minimizing the measurability loss is equivalent to minimizing the expected mean squared error. After discretization, \(y_0^\theta\) used by the measurability loss is estimated by
\begin{equation}
y_0^\theta \approx \phi(X^{(i)}_{t_H}) + \sum_{j=0}^{H-1}(Q(t_j,X^{(i)}_{t_j})+ \frac{1}{2}u^{(i)}(t_j,X^{(i)}_{t_j})^\intercal Ru^{(i)}(t_j,X^{(i)}_{t_j}))\Delta t - \sum_{j=k}^{H-1}\langle z^\theta(t_j,X^{(i)}_{t_j}),W_{t_{j+1}}-W_{t_j}\rangle.
\label{eq:discrete-def-ytheta}
\end{equation}

At last, we construct a new policy by minimizing the mean squared loss between it and \(-R^{-1}\Upsilon^\intercal z^\theta(t, x)\) along the trajectory \((t,X^{(i)}_t)\). That is, $u^{(i+1)}\equiv u^{\psi'}$, where $\psi'$ is obtained by
\begin{equation}
  \min_{\psi'\in\Psi} \operatorname{\mathbb{E}}\|u^{\psi'}(t,X^{(i)}_t) +R^{-1}\Upsilon^\intercal z^\theta(t,X_t^{(i)})\|^2.
\end{equation}

\begin{algorithm}
\DontPrintSemicolon
  \caption{Search the optimal state-feedback control policy \(u^*\)}

initialize neural networks $z^\theta$ and $u^\psi$\;
initialize zero policy $u^{(0)}(t,x) \equiv 0$\;
\uIf{model-based}{
set $\sigma, \Upsilon$ by Eq. \eqref{eq:model-based-settings}
}
\ElseIf{model-free}{
set $\sigma, \Upsilon$ by Eq. \eqref{eq:model-free-settings}
}
\ForEach{iteration $i$}{
\emph{\# solve forward SDE}\;
initialize a data buffer $\mathcal{D}$\;
\uIf{model-based}{
Generate sufficient trajectories $\{X^{(i)}_t,W^{(i)}_t\}$ via Eq. \eqref{eq:model-based-samples} with the current policy $u^{(i)}$.
}
\ElseIf{model-free}{
Generate sufficient trajectories $\{X^{(i)}_t,W^{(i)}_t\}$ via Eq. \eqref{eq:model-free-samples} with the current policy $u^{(i)}$.
}
store all generated trajectories in $\mathcal{D}$\;
\emph{\# policy evaluation}\;
\Repeat{$\theta$ convergence}{
random collect a batch of trajectories $\{X_t,W_t\}$ from buffer $\mathcal{D}$\;
calculate the measurability loss of this batch by Eq. \eqref{eq:discrete-def-ytheta}:
$$\mathit{MeasurLoss}(\theta)= \operatorname{Var} y_0^\theta$$
perform one gradient step to update $\theta$
}
\emph{\# policy improvement}\;
\Repeat{$\psi$ convergence}{
random collect a batch of trajectories $\{X_t,W_t\}$ from buffer $\mathcal{D}$\;
generate target control output $\{\hat{U}_t\}$ via Eq. \eqref{eq:policy-improvement-model-free}
$$\hat{U}_t =  -R^{-1}\Upsilon^\intercal z^\theta(t,X_t)$$
calculate the mean squared loss of this batch
$$\mathit{MeanSquare}(\psi) = \operatorname{\mathbb{E}}\|u^\psi(t,X_t) - \hat{U}_t\|^2$$
perform one gradient step to update $\psi$
}
set $u^{(i+1)} \leftarrow u^\psi$
}
\end{algorithm}

\textbf{Example 1.}  Consider a pure policy evaluation problem. Let the state process \(X_t\) be simply the \(n\)-dimensional standard Brownian motion, and running cost be equal to \(-n\) and terminal cost \(\phi(x) = x^\intercal x\). Now, one needs to determine the value function \(v(t,x)\). This is equivalent to consider the following BSDE
\begin{equation}
dY_t = n\,dt + \langle Z_t, dW_t\rangle,\qquad Y_T=\|W_T\|^2.
\label{eq:example-1}
\end{equation}
Assume that the vlaue function can be represented as \(\theta x^\intercal x\) with parameter \(\theta\in\mathbb{R}\).  One can validate this assumption by working out the conditional expectation in the definition of value function. Actually, \(v(t,x) = x^\intercal x\). Therefore, a rational loss function should have a minimizer \(\theta^* = 1\). To find out the optimal \(\theta^*\) in function class \(\{\theta x^\intercal x:\theta\in\mathbb{R}\}\), a common choice is to minimize the Monte Carlo loss, which is shown to be equivalent to the martingale loss proposed in \cite{Jia2021SSRN}. %
On the other hand, we choose to minimize the proposed measurability loss with \(z^\theta(t,x) = 2\theta x\), the derivative function class of \(\{\theta x^\intercal x: \theta\in\mathbb{R}\}\). At last, we try the Deep BSDE loss with the same choice of \(z^\theta\) plus a variable \(y_0^{DB}\) to be optimized simultaneously. We conduct this numerical experiment on PyTorch with a i7-7700 CPU. Initially, \(\theta\) is set to 0.5 and \(y_0^{DB}\) is 1.0. We simulate on time interval [0, 0.5 s] with time step \(\Delta t=0.01\). The optimizer is Adam implemented with learning rate 0.01. Each gradient step utilizes 32 state trajectory samples.

Figure. 1 plots the losses of above mentioned three approaches and the values of variables that are trained for different dimensions, \(n=1, 10, 100\). It can be seen that our approach basically coincides with the Deep BSDE approach, and that the martingale approach also converges but with a larger variance of loss function estimation. We observe that the martingale loss curve becomes more smooth if each gradient step uses 128 state trajectory samples. All the three approaches scale well with the state dimension \(n\), though the estimation of the initial value \(y_0^{DB}\) given by Deep BSDE in \(n=100\) is not good as the low-dimensional case.

\begin{figure}[htbp]
\centering
\includegraphics[width=.9\linewidth]{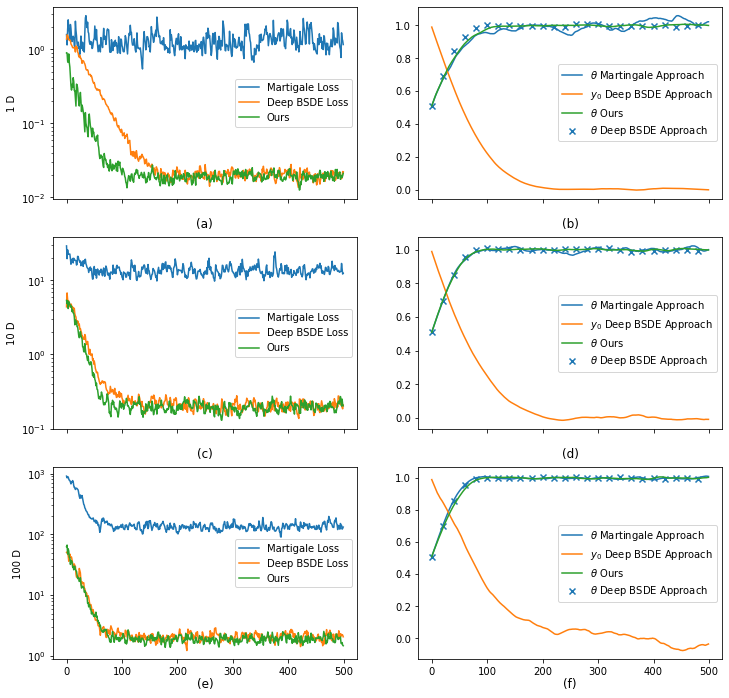}
\caption{Three different approaches for the linear BSDE with different dimensions. From the top to the bottom, the dimension \(n\) is 1, 10, 100, respectively. The three subplots in the left column show the validation loss during training and the three subplots in the right column show their estimation about \(\theta\). Deep BSDE method also estimate \(v(0,0)\) with its \(y_0\). The horizontal axis is training steps, which utilizes 32 trajectories to estimate the loss functions.}
\end{figure}

In practical applications, the form of true value function might be unknown. It is often the case that the function class \(\{(y^\theta, z^\theta):\theta\in\Theta\}\) does not contain the true answer \((v(t,x), \sigma^\intercal v_x(t,x))\). Nevertheless, we can define the optimal \(\theta^*\) by minimizing the following criteria
\begin{align}
\mathit{Yerr}(\theta) &= \operatorname{\mathbb{E}}\int_0^T|Y_t - y^\theta(t,X_t)|^2\,dt,\label{eq:yerr}\\
\mathit{Zerr}(\theta) &= \operatorname{\mathbb{E}}\int_0^T\|Z_t - z^\theta(t,X_t)\|^2\,dt,\label{eq:zerr}
\end{align}
where \((Y_t, Z_t)\) is the true solution of the BSDE. It is worth noting that the minimizers of these two criteria are sensitive to the choice of parameterizing strategy and do not necessarily coincide. Indeed, in this concrete example, if we choose \(\Theta=\mathbb{R}, y^\theta(t,x) = \theta \|x\|^4, z^\theta = 4\theta x\|x\|^2\) to solve BSDE \eqref{eq:example-1}, then it is clear that \((x^\intercal x, 2x)\) does not belong to this function class. After substituting the true answer \(Y_t=\|X_t\|^2, Z_t=2X_t\) and trial solutions \(y^\theta(\cdot,\cdot), z^\theta(\cdot,\cdot)\), we can obtain the analytic expression of the two criteria
\begin{align*}
\mathit{Yerr}(\theta) &= n(n+2)T^3\cdot\left[ \frac{1}{3} + \theta^2\cdot (n+4)(n+6)\cdot \frac{1}{5}T^2 - 2\theta\cdot (n+4)\cdot \frac{1}{4}T \right],\\
\mathit{Zerr}(\theta) &= 4nT^2\left[ \frac{1}{2} + \theta^2(n+2)(n+4)T^2 - \frac{4}{3}\theta(n+2)T \right].
\end{align*}
Therefore, we observe that the optimal \(\theta^*\) for fitting value function is \(\theta^*_Y = \frac{5}{4(n+6)T}\), whereas the optimal \(\theta^*\) for fitting value function gradient is \(\theta^*_Z = \frac{2}{3(n+4)T}\). We rerun the program in aforementioned numerical experiment except for \(y^\theta(t,x)=\theta \|x\|^4, z^\theta(t,x)=4x\|x\|^2\) this time. The numerical results are given in Figure 2.

Compared with Figure 1, Figure 2 shows that the martingale approach does not converge to the same solution as that of the other two methods. We plot the theoretical value \(\theta^*_Y, \theta^*_Z\) to stress this. It it clear that the martingale approach actually converges to \(\theta^*_Y\) while the Deep BSDE method and ours both converge to \(\theta^*_Z\). This is expected because Theorem 1 ensures that \(\mathit{Zerr}\) is strictly equal to the measurability loss we proposed, and Theorem 2 ensures that minimizing this martingale loss is equivalent to the Deep BSDE method. Moreover, it has also been proved in \cite{Jia2021SSRN} that \(\mathit{Yerr}\) is strictly equal to the loss used in martingale approach plus a constant term. This numerical experiment suggests that the martingale approach and our method are two parallel ways for solving linear BSDE, one for \(Y_t\) and the other for \(Z_t\).

\begin{figure}[htbp]
\centering
\includegraphics[width=.9\linewidth]{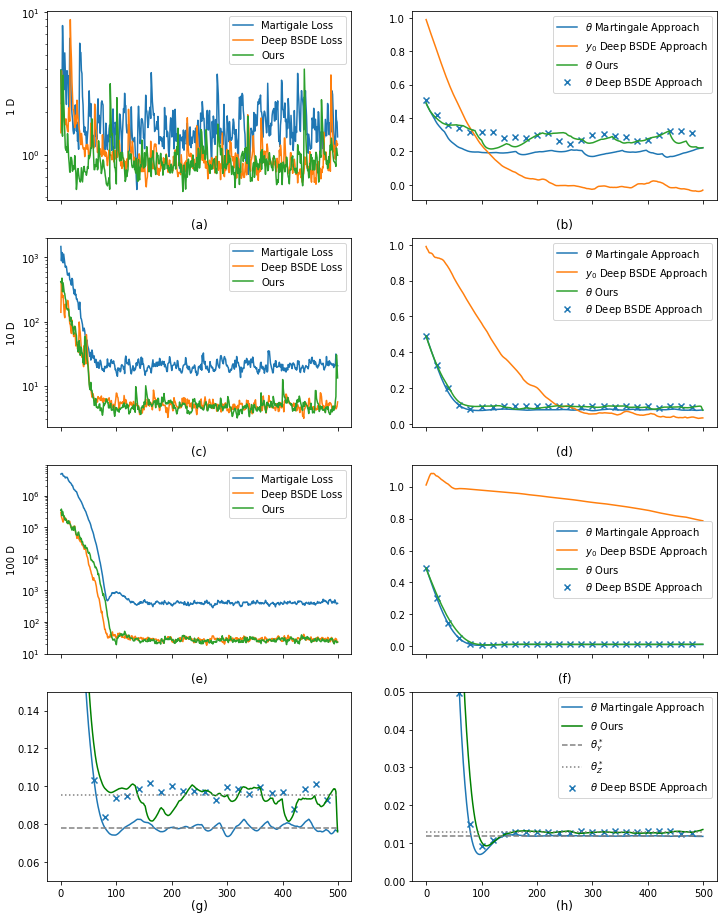}
\caption{Rerun with different parameterizing strategy. The first three rows have the same meaning as those of Figure 1. To compare with the theoretical prediction, we rescale (d), (e) along vertical axis and add two horizontal line for \(\theta^*_Y, \theta^*_Z\) in subplots (g), (h) respectively.}
\end{figure}

\textbf{Example 2}. Consider the task of swinging up a pendulum studied in \cite{Lee2021,Pereira2019,ISOCPI2017}. This system suits the formulation of Eq. \eqref{eq:dyn} with state dimension \(D_x=2\) and control dimension \(D_u=1\). The model-specific term are (let \((\theta_p,\dot{\theta}_p)^\intercal\) denote the system state)
\begin{align*}
F(t,(\theta_p,\dot{\theta}_p)^\intercal) &= (\dot{\theta}_p, (a\sin\theta_p - b\dot{\theta}_p)/I_{\mathrm{inertia}})^\intercal,\\
G(t,(\theta_p,\dot{\theta}_p)^\intercal) &= (0, \cos\theta_p / I_{\mathrm{inertia}})^\intercal.
\end{align*}
For both the model-based and model-free settings, we choose \(\sigma_0=1.414\) and the diffusion coefficient \(\sigma(t,x)\) is determined by Eq. \eqref{eq:model-based-settings} and \eqref{eq:model-free-settings} accordingly. We complete the description of forward dynamics Eq. \eqref{eq:sto-dyn} by setting \(a=9.8,b=0.1,I_{\mathrm{inertia}}=1.0\). For the cost functional Eq. \eqref{eq:op}, swinging up task requires to drive the system from initial state \(x_0=(\pi,0)^\intercal\) to origin \(x^*=(0,0)^\intercal\) with minimal time and energy. Thus we choose \(\phi\equiv 0, Q(t,x) = (x-x^*)^\intercal\Lambda(x-x^*), R=0.005I\), where \(\Lambda\) is a diagonal matrix with coefficient \(\lambda_1=1.01, \lambda_2=0.01\). The time interval is \([0,1.0\text{ s}]\) and time step is \(\Delta t=0.01\text{ s}\). When implement the proposed policy iteration algorithm, the size of buffer \(\mathcal{D}\) is 12800, and batch size is 128. \(z^\theta(\cdot,\cdot)\) and \(u^\psi(\cdot,\cdot)\) are represented by feedforward neural networks with a single hidden layer plus a batch norm layer. The hidden size is 16 and activation function is \(\operatorname{tanh}\). The whole program is built on PyTorch and the optimizer is Adam with learning rate \(1\times10^{-4}\) and weight decay \(1\times10^{-8}\).

Figure 3 plots the obtained controlled trajectory of 4 successive policy iterations on the original deterministic system. It shows that our policy iteration algorithm works well with or without system models. It can be seen that 2 or 3 policy iteration leads to a considerable good policy. This is consistent with the experience in discrete time reinforcement learning problems, where it has been proved that worst case of policy iteration approach is slow but in practice it needs only very few iterations \cite{Bertsekas2019}.

\begin{figure}[htbp]
\centering
\includegraphics[width=.9\linewidth]{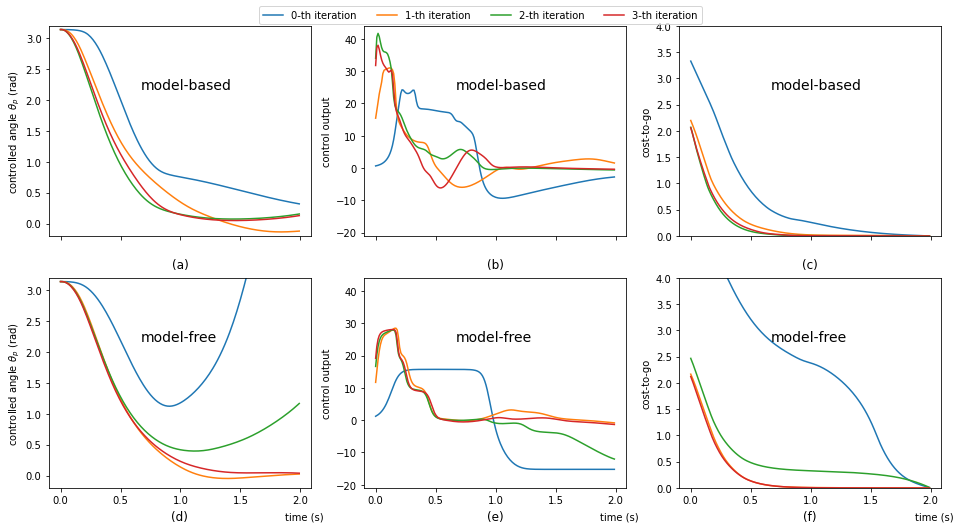}
\caption{Performance in 4 successive iterations on the original deterministic system. The left, middle and right columns are controlled angle \(\theta_p\), control output and cost-to-go of the obtained trajectories. The first row and second row are plotted with data generated by model-based and model-free algorithm respectively.}
\end{figure}

\section{Conclusion}
\label{sec:orgde49b96}

In this paper, we study a class of nonlinear optimal control problems through the lens of PDEs and FBSDEs. Inspired from reinforcement learning, we start by introducing some exploration noise and transform the original problem to the stochastic setting. Then, a policy iteration framework is developed for diffusion process that is parallel to that for MDP. By exploiting the quadratic structure of Hamiltonian, the value function gradients are essential to policy improvement instead of the value function itself. To evaluate the policy, a modified Deep BSDE method with  measurability loss is introduced, and this loss is proved to be strictly equal to the expected mean squared error of the diffusion term and its approximation. Under this framework, both a model-based and a model-free numerical schemes  are proposed to search the optimal control. Numerical experiment results on a simple BSDE are in good agreement on our theory of Deep BSDE-ML method, and the nonlinear optimal control example demonstrates that the model-free approach performs as good as its model-based counterpart. For future research, a theoretically interesting one is how to extend the considered problem to the infinite horizon and to develop online learning algorithms.

\bibliographystyle{plain}
\bibliography{dbpg-pi}

\section*{Appendix}
\label{sec:orged7ce40}

\textbf{Proof of Theorem 1}. By Itô's lemma, the initial value can be rewritten as
\begin{equation*}
v^\psi(0,x_0) = v^\psi(T,X_T) + \int_0^Tg^\psi(t,X_t)\,dt - \int_0^T\langle \sigma^\intercal v_x^\psi(t,X_t),dW_t\rangle,
\end{equation*}
where \(\{X_t, 0\leq t\leq T\}\) is the solution of the forward dynamics driven by \(u^\psi\). On the other hand, by the conditional expectation definition of value function, we have
\begin{align*}
\operatorname{\mathbb{E}}y_0^\theta &= \operatorname{\mathbb{E}}\left[\phi(X_T) + \int_0^Tg^\psi(t,X_t)\,dt - \int_0^T\langle z^\theta(t,X_t),dW_t\rangle\right]\\
&= \operatorname{\mathbb{E}}\left[\phi(X_T) + \int_0^Tg^\psi(t,X_t)\,dt\rangle\right] \\
&= \operatorname{\mathbb{E}}\left[\phi(X_T) + \int_0^Tg^\psi(t,X_t)\,dt\rangle\biggm| X_0\right]\\
&= v^\psi(0,x_0).
\end{align*}
The thrid equality comes from the fact that the \(\sigma\)-algebra generated by \(X_0\) is trivial because \(X_0\) is deterministic and equal to \(x_0\). The second equality holds because \(\int_0^s\langle z^\theta(t,X_t),dW_t\rangle\) is a martingale. Expanding the left hand side of Eq. \eqref{eq:z-err} yields the desired expression
\begin{align*}
\operatorname{\mathbb{E}}|y_0^\theta - \operatorname{\mathbb{E}}y_0^\theta|^2 &= \operatorname{\mathbb{E}}|y_0^\theta - v^\psi(0,x_0)|^2\\
&= \operatorname{\mathbb{E}}\biggl|\phi(X_T) + \int_0^Tg^\psi(t,X_t)\,dt - \int_0^Tz^\theta(t,X_t),dW_t\rangle - v^\psi(T,X_T)\\&\qquad - \int_0^Tg^\psi(t,X_t)\,dt + \int_0^T\langle\sigma^\intercal v_x^\psi(t,X_t),dW_t\rangle \biggr|^2\\
&= \operatorname{\mathbb{E}}\biggl|\int_0^T\langle z^\theta(t,X_t) - \sigma^\intercal v_x^\psi(t,X_t),dW_t\rangle\biggr|^2\\
&= \operatorname{\mathbb{E}}\int_0^T\|z^\theta(t,X_t) - \sigma^\intercal v_x^\psi(t,X_t)\|^2\,dt \qquad\text{(Itô's Isometry)}.\tag*{\(\square\)}
\end{align*}

\textbf{Proof of Theorem 2}. It is sufficient to prove the following equality holds for any \(\theta\) and \(y_0^{DB}\),
\begin{equation}
\operatorname{\mathbb{E}}|y_T^{\theta, DB} - \phi(X_T)|^2 = |y_0^{DB} - v^\psi(0,x_0)|^2 + \operatorname{\mathbb{E}}|y_0^\theta - \operatorname{\mathbb{E}}y_0^\theta|^2.
\label{eq:app-tmp-1}
\end{equation}
Eq. \eqref{eq:app-tmp-1} says that the Deep BSDE loss Eq. \eqref{eq:def-DB-loss} can be decomposed into the proposed measurability loss Eq. \eqref{eq:def-measurability-loss} plus \(|y_0^{DB} - v^\psi(0,x_0)|^2\), which is independent of \(\theta\) and attains 0 at the unique minimizer \(y_0^{*,DB}=v^\psi(0,x_0)\). Therefore, once the Deep BSDE loss attains its minimum at \((y_0^{*,DB}, \theta^*)\), it must hold that \(y_0^{*,DB}=v^\psi(0,x_0)\) and \(\theta^*\) is a minimizer of the measurability loss. Convsersely, if there exists a \(\theta^*\) minimizing the measurability loss, then \((v^\psi(0,x_0), \theta^*)\) minimizes the Deep BSDE loss.

Now we prove Eq. \eqref{eq:app-tmp-1}. Expanding the definition of the auxiliary process \(y_t^{\theta,DB}\) yields
\begin{align*}
\operatorname{\mathbb{E}}|y_T^{\theta,DB} - \phi(X_T)|^2 &= \operatorname{\mathbb{E}}\biggl|y_0^{DB} -\int_0^Tg^\psi(t,X_t)\,dt + \int_0^T\langle z^\theta(t,X_t),dW_t\rangle - \phi(X_T)\biggr|^2\\
&= \operatorname{\mathbb{E}}\biggl|\phi(X_T) + \int_0^Tg^\psi(t,X_t)\,dt - \int_0^T\langle z^\theta(t,X_t),dW_t\rangle - y_0^{DB}\biggr|^2.
\end{align*}
Then we substitute the definition of the auxiliary process \(y_t^{\theta}\)
\begin{align}
\operatorname{\mathbb{E}}|y_T^{\theta,DB} - \phi(X_T)|^2 &= \operatorname{\mathbb{E}}|y_0^{\theta} - y_0^{DB}|^2 \notag\\
&= \operatorname{\mathbb{E}}|y_0^{\theta} - \operatorname{\mathbb{E}}y_0^\theta + \operatorname{\mathbb{E}}y_0^\theta - y_0^{DB}|^2 \notag\\
&= \operatorname{\mathbb{E}}|y_0^{\theta} - \operatorname{\mathbb{E}}y_0^\theta|^2 + | \operatorname{\mathbb{E}}y_0^\theta - y_0^{DB}|^2 + 2\bigl(\operatorname{\mathbb{E}}y_0^\theta - y_0^{DB}\bigr)\cdot\operatorname{\mathbb{E}}\bigl[y_0^\theta - \operatorname{\mathbb{E}}y_0^\theta\bigr] \notag\\
&= |y_0^{DB} - \operatorname{\mathbb{E}}y_0^\theta|^2 + \operatorname{\mathbb{E}}|y_0^\theta - \operatorname{\mathbb{E}}y_0^\theta|^2.
\label{eq:app-tmp-2}
\end{align}
Following the same argument in proof of Theorem 1, we have \(\operatorname{\mathbb{E}}y_0^\theta = v^\psi(0,x_0)\). Replacing \(\operatorname{\mathbb{E}}y_0^\theta\)  in Eq. \eqref{eq:app-tmp-2} with \(v^\psi(0,x_0)\) results in Eq. \eqref{eq:app-tmp-1}. This completes the proof. \hfill\(\square\)

\textbf{Proof of Theorem 3}. By applying Itô's formula, we obtain
\begin{align}
dv^\psi(t,X^\psi_t) &= \mathcal{L}^\psi v^\psi(t,X^\psi_t)\,dt + \langle \sigma^\intercal v^\psi_x(t,X_t^\psi),dW_t\rangle,\label{eq:app-tmp-3}\\
dv^\theta(t,X^\psi_t)&= \mathcal{L}^\psi v^\theta(t,X_t^\psi)\,dt + \langle \sigma^\intercal v^\theta_x(t,X_t^\psi),dW_t\rangle.\label{eq:app-tmp-4}
\end{align}
Subtracting Eq. \eqref{eq:app-tmp-3} from Eq. \eqref{eq:app-tmp-4} and  integrating on \([s,T]\) results in
\begin{align*}
v^\psi(s,X_s^\psi) - v^\theta(s,X_s^\psi) &= v^\psi(T,X^\psi_T) - \int_s^T\mathcal{L}^\psi v^\psi(t,X_t^\psi)\,dt - \int_s^T\langle \sigma^\intercal v^\psi_x(t,X_t^\psi),dW_t\rangle \\
&\qquad -v^\theta(T,X^\psi_T) + \int_s^T\mathcal{L}^\psi v^\theta(t,X_t^\psi)\,dt + \int_s^T\langle \sigma^\intercal v^\theta_x(t,X_t^\psi),dW_t\rangle \\
&= \int_s^T\mathcal{L}^\psi(v^\theta - v^\psi)(t,X_t^\psi)\,dt + \int_s^T\langle \sigma^\intercal(v_x^\theta - v_x^\psi)(t,X_t^\psi),dW_t\rangle.
\end{align*}
Taking expectation on both sides and noting that \(X_s^\psi\) is deterministic and equal to \(x_s\), there is
\begin{equation}
v^\psi(s,x_s) - v^\theta(s,x_s) = \operatorname{\mathbb{E}}\int_s^T\mathcal{L}^\psi(v^\theta - v^\psi)(t,X_t^\psi)\,dt.
\label{eq:app-thm3-tmp2}
\end{equation}
The second term disappears because \(\int_s^\tau\langle \sigma^\intercal (v_x^\theta - v_x^\psi)(t,X_t^\psi),dW_t\rangle\) is a martingale. We complete the proof by showing that
\begin{equation}
\mathcal{L}^\psi(v^\theta - v^\psi)(t,x) = H(t,x,u^\psi(t,x),v_x^\theta(t,x)) - H(t,x,u^\theta(t,x),v_x^\theta(t,x)),\qquad \forall t,x.
\label{eq:app-thm3-tmp1}
\end{equation}
By exploiting the linear dependency on control of the system dynamics and the fact that the diffusion coefficient is independent of control, the infinitesimal generator \(\mathcal{L}^\psi\) can be represented with \(\mathcal{L}^\theta\) plus a compensate term
\begin{align*}
\mathcal{L}^\psi v^\theta(t,x) &= v_t^\theta(t,x) + \langle v_x^\theta(t,x), F(t,x) + G(t,x)u^\psi(t,x)\rangle + \frac{1}{2}\operatorname{tr}\{\sigma\sigma^\intercal v^\theta_{xx}(t,x)\}\\
&= v_t^\theta(t,x) + \langle v_x^\theta(t,x), F(t,x) + G(t,x)u^\theta(t,x)\rangle + \frac{1}{2}\operatorname{tr}\{\sigma\sigma^\intercal v^\theta_{xx}(t,x)\}\\
&\qquad\qquad + \langle v_x^\theta(t,x), G(t,x)u^\psi(t,x) - G(t,x)u^\theta(t,x)\rangle \\
&= \mathcal{L}^\theta v^\theta(t,x) + \langle v_x^\theta(t,x),G(t,x)u^\psi(t,x) - G(t,x)u^\theta(t,x)\rangle.
\end{align*}
Recall that the PDE formulation of value function means that  \(\mathcal{L}^\theta v^\theta + g^\theta = 0\) and \(\mathcal{L}^\psi v^\psi + g^\psi = 0\), the left hand side of Eq. \eqref{eq:app-thm3-tmp1} can be rewritten as
\begin{align*}
\mathcal{L}^\psi(v^\theta - v^\psi)(t,x) &= \mathcal{L}^\theta v^\theta (t,x) - \mathcal{L}^\psi v^\psi(t,x) + \langle v_x^\theta(t,x), G(t,x)u^\psi(t,x) - u^\theta(t,x)\rangle \\
&= -g^\theta(t,x) + g^\psi(t,x) + \langle v_x^\theta(t,x), F(t,x) + G(t,x)u^\psi(t,x)\rangle\notag\\
&\qquad - \langle v_x^\theta(t,x),F(t,x) + G(t,x)u^\theta(t,x)\rangle\\
&=H(t,x,u^\psi(t,x),v_x^\theta(t,x)) - H(t,x,u^\theta(t,x),v_x^\theta(t,x)).
\end{align*}
Thus Eq. \eqref{eq:app-thm3-tmp1} holds and substituting it into Eq. \eqref{eq:app-thm3-tmp2} yields the conclusion.\hfill\(\square\)

\end{document}